\documentclass[12pt]{amsart}
    \title[Generalised Moore Spectra]{Generalised Moore Spectra in a Triangulated Category}
    \author{David Pauksztello}
    \date{2 June 2010}
    
\usepackage{amssymb}
\usepackage{latexsym}
\usepackage{amsmath}
\usepackage{mathrsfs}
\usepackage[all]{xy}

\newtheorem{dfn}{Definition}[section]
\newtheorem{thrm}[dfn]{Theorem}

\newtheorem{lem}[dfn]{Lemma}
\newtheorem{prp}[dfn]{Proposition}
\newtheorem{corr}[dfn]{Corollary}

\newtheorem{rmk}[dfn]{Remark}
\newtheorem{rmks}[dfn]{Remarks}

\newtheorem{set}[dfn]{Setup}
\newtheorem{example}[dfn]{Example}
\newtheorem{examples}[dfn]{Examples}
\newtheorem{observe}[dfn]{Observation}
\newtheorem{definitions}[dfn]{Definitions}
\newtheorem{conject}[dfn]{Conjecture}

\newenvironment{prf}{\noindent\textbf{Proof: }}{\hfill $\Box$\ }
\newenvironment{enumerateroman}{\begin{enumerate}}{\end{enumerate}}

\newcommand{\df}[1]{\begin{dfn}\emph{#1}}
\newcommand{\edf}{\end{dfn}}
\newcommand{\thm}{\begin{thrm}}
\newcommand{\ethm}{\end{thrm}}
\newcommand{\lemma}{\begin{lem}}
\newcommand{\elemma}{\end{lem}}
\newcommand{\prop}{\begin{prp}}
\newcommand{\eprop}{\end{prp}}
\newcommand{\rem}[1]{\begin{rmk}\emph{#1}}
\newcommand{\erem}{\end{rmk}}
\newcommand{\rems}[1]{\begin{rmks}\emph{#1}}
\newcommand{\erems}{\end{rmks}}
\newcommand{\pf}{\begin{prf}}
\newcommand{\epf}{\end{prf}}
\newcommand{\cor}{\begin{corr}}
\newcommand{\ecor}{\end{corr}}
\newcommand{\setup}[1]{\begin{set}\emph{#1}}
\newcommand{\esetup}{\end{set}}
\newcommand{\eg}[1]{\begin{example}\emph{#1}}
\newcommand{\eeg}{\end{example}}
\newcommand{\egs}[1]{\begin{examples}\emph{#1}}
\newcommand{\eegs}{\end{examples}}
\newcommand{\obs}[1]{\begin{observe}\emph{#1}}
\newcommand{\eobs}{\end{observe}}
\newcommand{\dfs}[1]{\begin{definitions}\emph{#1}}
\newcommand{\edfs}{\end{definitions}}
\newcommand{\conj}{\begin{conject}}
\newcommand{\econj}{\end{conject}}

\renewcommand{\geq}{\geqslant}
\renewcommand{\leq}{\leqs}
\renewcommand{\phi}{\varphi}

\DeclareMathAlphabet{\mathpzc}{OT1}{pzc}{m}{it}
\newcommand{\D}{\mathscr{D}}

\newcommand{\T}{\mathscr{T}}

\newcommand{\A}{\mathscr{A}}
\newcommand{\B}{\mathscr{B}}
\newcommand{\C}{\mathscr{C}}

\newcommand{\M}{\mathscr{M}}
\newcommand{\N}{\mathscr{N}}
\newcommand{\Abel}{\mathsf{Ab}}

\newcommand{\Hom}[3]{\mathsf{Hom}_{#1}(#2,#3)}

\newcommand{\rightiso}{\stackrel{\sim}{\longrightarrow}}
\newcommand{\Ext}[2]{\mathsf{Ext}^{#1}_{#2}}

\newcommand{\Endo}[1]{\mathsf{End}(#1)}
\newcommand{\Mod}[1]{\mathsf{Mod}(#1)}
\renewcommand{\mod}[1]{\mathsf{mod}(#1)}

\newcommand{\Proj}[1]{\mathsf{Proj}(#1)}
\newcommand{\proj}[1]{\mathsf{proj}(#1)}
\newcommand{\Add}[1]{\mathsf{Add}(#1)}
\newcommand{\Projk}[2]{\mathsf{Proj}^{#1}(#2)}
\newcommand{\projk}[2]{\mathsf{proj}^{#1}(#2)}
\newcommand{\op}{\mathsf{op}}
\newcommand{\Sop}{S^{\op}}
\newcommand{\EndC}{\mathsf{End}_{\T}(C)}
\newcommand{\rightlabel}[1]{\stackrel{#1}{\longrightarrow}}

\newcommand{\Essential}[1]{\mathsf{Ess.Im}\, #1}

\newcommand{\Z}{\mathbb{Z}}

\newcommand{\gldim}{\mathsf{gldim}\,}
\newcommand{\lgldim}{\mathsf{l.gldim}\,}
\newcommand{\rgldim}{\mathsf{r.gldim}\,}
\newcommand{\projdim}[1]{\mathsf{projdim}_{#1}\,}
\newcommand{\tri}[3]{#1\rightarrow #2\rightarrow #3\rightarrow \Sigma #1}
\newcommand{\trilabel}[4]{#1\stackrel{#4}{\longrightarrow} #2\longrightarrow #3\longrightarrow \Sigma #1}
\newcommand{\trilabels}[6]{#1\stackrel{#4}{\longrightarrow} #2\stackrel{#5}{\longrightarrow} #3\stackrel{#6}{\longrightarrow} \Sigma #1}

\newcommand{\leqs}{\leqslant}

\entrymodifiers={+!!<0pt,\fontdimen22\textfont2>}
   
\numberwithin{equation}{section}

\begin{document}

\begin{abstract}
In this paper we consider a construction in an arbitrary triangulated category $\T$ which resembles the notion of a Moore spectrum in algebraic topology. Namely, given a compact object $C$ of $\T$ satisfying some finite tilting assumptions, we obtain a functor which ``approximates'' objects from the module category of the endomorphism algebra of $C$ in $\T$. This generalises and extends a construction of J\o rgensen which appears in \cite{Jorgensen2} in connection with lifts of certain homological functors of derived categories. We show that this new functor is well-behaved with respect to short exact sequences and distinguished triangles, and as a consequence we obtain a new way of embedding a module category in a triangulated category. As an example of the theory, we recover Keller's canonical embedding of the module category of a path algebra of a quiver with no oriented cycles into its $u$-cluster category of $u\geqslant 2$.
\end{abstract}
	
\address{Institut f\"ur Algebra, Zahlentheorie und Diskrete Mathematik, Fakult\"at f\"ur Mathematik und Physik, Leibniz Universit\"at Hannover, Welfengarten 1, 30167 Hannover, Germany.}
\email{pauk@math.uni-hannover.de}
    
    
\maketitle

\section*{Introduction}

In this paper we discuss the existence of ``Moore spectra'' in a triangulated category. The terminology ``Moore spectra'' employed in this paper is borrowed from algebraic topology, see \cite{Margolis}. While the notion discussed here does not coincide with its counterpart in algebraic topology its spirit is the same.

In algebraic topology, the notion of spectra can be considered as one of ``generalised topological spaces''. In this setting one uses the idea of a Moore spectrum to construct a spectrum with a single (pre-defined) non-vanishing homology group; c.f. the notion of an Eilenberg-MacLane space for homotopy groups.  For instance, suppose $A$ is an abelian group, the Moore spectrum $MA$ of $A$ is a spectrum with
$$
H(\Sigma^{i}MA)=
\left\{
\begin{array}{ll}
 A & i=0 \\
0 & i\neq 0,
\end{array}
\right.
$$
where $\Sigma$ is the suspension functor in the category of spectra.

Analogously, in this paper we shall consider the following setup: suppose $C$ is a compact object of a triangulated category $\T$ which has set indexed coproducts satisfying some finite tilting assumptions (see Setup \ref{mainsetup} for precise conditions). Then consider the endomorphism algebra $\Sop=\Endo{C^{\op}}$ of $C$. We are looking for ways to ``represent'' or ``approximate'' an $\Sop$-module in the triangulated category $\T$. Compare this with the idea of Moore spectrum in algebraic topology where one starts with given homology groups and ``approximates'' a representative topological spae with given homology.

The idea of ``Moore spectra'' in triangulated categories in general was first studied by J\o rgensen in \cite{Jorgensen2}; J\o rgensen's construction is also in the same spirit as its namesake in algebraic topology. It is used as a tool when considering certain lifts of homological functors between the derived category of the integers $\D(\Z)$ and an arbitrary triangulated category $\T$.

However, the notion of Moore spectra in a triangulated category as developed by J\o rgensen in \cite{Jorgensen2} holds independent interest, in particular it yields a new technique of embedding an abelian category in a triangulated category and of obtaining a module category from a triangulated category in a nice way. The construction given in this paper is a higher analogue of J\o rgensen's construction, though the language and techniques used to prove the theorems are very different.

The outline of this paper is as follows: Section \ref{notation} recalls some definitions and fixes notation; Section \ref{jorgensen_construction} consists of a brief review of J\o rgensen's construction as it appears in \cite{Jorgensen2} and in Section \ref{general_construction} we present the general construction and state the main theorem (Theorem \ref{moorespectra}) of this paper. Section \ref{proof} is devoted to the proof of the main theorem. In Section \ref{revisited} we return briefly to J\o rgensen's construction and in Section \ref{well_behaved} we show that the functor obtained in Theorem \ref{moorespectra} is well behaved, in particular, that it takes short exact sequences to distinguished triangles and induces a natural Hom-Ext sequence. The final section, Section \ref{cluster_categories} contains an example: we show how the Moore spectra functor recovers Keller's canonical embedding of the module category of a hereditary algebra into its $u$-cluster category for $u\geqslant 2$.

Definitions and examples of triangulated categories can be found in \cite{Hartshorne} and \cite{Neeman2}. Background in homological algebra can be found in \cite{Hilton} and \cite{Weibel}, and in representation theory and algebra in \cite{Assem}.

\section{Notation and preliminaries}\label{notation}

Throughout this paper, unless stated otherwise, $\T$ will be a triangulated category with set indexed coproducts.

\subsection{Additive closure and compact objects}

By $\Add{C}$ we denote the \textit{(infinite) additive closure} of $C$ in $\T$, that is the smallest full subcategory of $\T$ whose objects are direct summands of (possibly infinite) set indexed coproducts of $C$.

Recall the following definition from \cite{Neeman}:

\df
{An object $C$ of $\T$ is called \textit{compact} if the functor $\Hom{\T}{C}{-}$ commutes with set indexed coproducts; that is given a set of objects $\{X_{i}\}_{i\in I}$ in $\T$ then there is a canonical isomorphism
$$\coprod_{i\in I}\Hom{\T}{C}{X_{i}}\rightiso\Hom{\T}{C}{\coprod_{i\in I}X_{i}}.$$}
\edf

\subsection{Projective dimension and global dimension}

Let $R$ be a ring. We shall denote by $\Mod{R}$ the category of left $R$-modules. We shall denote the category of right $R$-modules by $\Mod{R^{\op}}$. We shall refer to a left $R$-module simply as an $R$-module and a right $R$-module simply as an $R^{\op}$-module.

Recall that an $R$-module $P$ is called \textit{projective} if the functor $\Hom{R}{P}{-}:\Mod{R}\rightarrow\Abel$ is exact, that is the sequence induced by the application of the functor on a short exact sequence is also a short exact sequence.

Let $A$ be an $R$-module. A \textit{projective resolution} of $A$ consists of chain complex
$$P:\, \cdots\rightarrow P_{n}\rightarrow P_{n-1}\rightarrow\cdots\rightarrow P_{1}\rightarrow P_{0}\rightarrow 0$$
such that $H_{i}(P)=0$ for $i\geqslant 1$ and together with an isomorphism $H_{0}(P)\rightiso A$, where $H_{i}(P)$ denotes the $i^{\textnormal{th}}$-homology group of the complex $P$.

\df
{Let $A$ be an $R$-module. The \textit{projective dimension} of $A$, written $\projdim{R}A$, is the smallest integer $n$ such that there is a projective resolution of $A$,
$$0\rightarrow P_{n}\rightarrow P_{n-1}\rightarrow\cdots\rightarrow P_{1}\rightarrow P_{0}\rightarrow A.$$
We write $\projdim{R}A=n$. If no finite resolution exists we say $A$ has infinite projective dimension.}
\edf

\df
{Let $R$ be a ring. The \textit{left global dimension} of $R$ is defined as
$$\lgldim R :=\sup\{\projdim{R}M\, |\, M\in\Mod{R}\}.$$}
\edf

The \textit{right global dimension} of $R$, denoted $\rgldim R$, is defined similarly. Where right or left global dimension can be understood from context we shall simply write $\gldim R$ and refer to the \textit{global dimension} of $R$.

\df
{We define the following full subcategories of $\Mod{R}$: 
\begin{itemize}
\item[\textit{(a)}] By $\Proj{R}$ we denote the full subcategory of $\Mod{R}$ whose objects are all projective left $R$-modules. 
\item[\textit{(b)}] By $\Projk{k}{R}$ we denote the full subcategory of $\Mod{R}$ whose objects are left $R$-modules with projective dimension at most $k$.
\end{itemize}}
\edf

\subsection{Adjoint functors}

Recall that for two categories $\A$ and $\B$ an \textit{adjunction} is a pair of functors $F$ and $G$,
$$\xymatrix{\A\ar@<-1ex>[rr]_-{G} & & \B,\ar@<-1ex>[ll]_-{F}}$$
such that there is a natural isomorphism
$$\Hom{\A}{FB}{A}\cong\Hom{\B}{B}{GA}.$$
The functor $F$ is called the \textit{left adjoint} and the functor $G$ is called the \textit{right adjoint}. $F$ (resp. $G$) is said to be \textit{left (resp. right) adjoint} to $G$ (resp. $F$).

Recall that a functor $F:\A\rightarrow \B$ is called \textit{fully faithful} (or \textit{full and faithful}) if
$$\Hom{\A}{A}{A'}\cong\Hom{\B}{FA}{FA'}$$
for all objects $A$ and $A'$ of $\A$.

We shall need the following observation from \cite{MacLane} connecting the notions of fully faithfulness with adjoint functors. 

\lemma
[\cite{MacLane}, Theorem IV.3.1]
Let $\A$ and $\B$ be categories and suppose
$$\xymatrix{\A\ar@<-1ex>[rr]_-{G} & & \B,\ar@<-1ex>[ll]_-{F}}$$
is an adjunction with $F$ left adjoint to $G$. Then the unit of the adjunction is an isomorphism if and only if $F$ is fully faithful.
\label{unit_isomorphism}
\elemma

We are now ready to review J\o rgensen's construction.

\section{J\o rgensen's construction}\label{jorgensen_construction}

Let $R$ be a commutative ring. Recall that a triangulated category $\T$ is called \textit{$R$-linear} if for any two objects $X$ and $Y$ of $\T$ the Hom-space $\Hom{\T}{X}{Y}$ is an $R$-module and the composition of morphisms is $R$-bilinear. A functor $F:\T\rightarrow\T'$ of $R$-linear triangulated categories is said to be $R$-linear if $F(r\alpha)=rF(\alpha)$ for all morphisms $\alpha$ of $\T$ and all elements $r\in R$;  see \cite[Definition 1.2]{Jorgensen2}.

In \cite[Section 4]{Jorgensen2}, J\o rgensen considers the following setup. 

\newtheorem{joder}[dfn]{Setup}
\begin{joder}[\cite{Jorgensen2}, Setups 4.1 and 4.12]
{Let $R$ be a principal ideal domain, let $\T$ be an $R$-linear triangulated category with set indexed coproducts and let $C$ be a compact object of $\T$ which satisfies the following assumptions:
\begin{enumerate}
\item $\Hom{\T}{C}{C}$ is a flat $R$-module;
\item $\Hom{\T}{C}{\Sigma^{-1}C}=0$;
\item $\Hom{\T}{C}{\Sigma^{-2}C}=0$.
\end{enumerate}}
\label{jorgensen_setup}
\end{joder}

The idea in \cite{Jorgensen2} is to construct the best possible approximation of an $R$-module $A$ in $\T$. This approximation of $A$ in $\T$ is denoted by $M(A)$ and is called the \textit{Moore spectrum} of $A$ in $\T$. J\o rgensen's construction depends on which object, $C$, of $\T$ is employed as the Moore spectrum of the ring itself.

J\o rgensen introduces an auxilliary category $\M$, an analogue of which we shall also introduce in the general construction. Below is the definition of $\M$ in J\o rgensen's setting; see \cite[Definition 4.3]{Jorgensen2}.

\df
{Let $\M$ be the full subcategory of $\T$ consisting of all objects $M$ of $\T$ which occur in distinguished triangles of the form
$$\trilabel{C\otimes F_{1}}{C\otimes F_{0}}{M}{1_{C}\otimes f},$$
when $A$ is an $R$-module with free resolution
$$0\longrightarrow F_{1}\rightlabel{f}F_{0}\longrightarrow A\longrightarrow 0.$$}
\label{jorgensen_M}
\edf

Note that in this definition, the tensor product is not the ``usual'' tensor product but is in fact a bifunctor $-\otimes-:\T\times\mathsf{Free}(R)\rightarrow \T$ which is $R$-linear, preserves set indexed coproducts and has $X\otimes R\cong X$ for each $X$ in $\T$. Here, $\mathsf{Free}(R)$ denotes the full subcategory of $\Mod{R}$ consisting of all free $R$-modules. Indeed, the construction of $-\otimes-$, for each $X$ in $\T$ and $F$ in $\mathsf{Free}(R)$, identifies $X\otimes F$ with a coproduct $\coprod_{I_{F}}X$, where $I_{F}$ is an indexing set for a basis of $F$. For full details of the construction see \cite[Construction 1.4 and Lemma 1.5]{Jorgensen2}.
J\o rgensen then obtains the following theorem.

\thm
[\cite{Jorgensen2}, Proposition 4.7 and Theorem 4.9]
Under the hypotheses of Setup \ref{jorgensen_setup}, the functor
$$\Hom{\T}{C}{-}:\M\rightarrow\Mod{R}$$
has a left adjoint
$$M:\Mod{R}\rightarrow\M.$$
If $M$ is viewed as a functor
$M:\Mod{R}\rightarrow \T$
by composition with the inclusion functor $i:\M\hookrightarrow\T$, then $M$ is an $R$-linear functor, it has $M(R)\cong C$ and it preserves set indexed coproducts.
\label{jorgensen_thm}
\ethm

J\o rgensen then continues to prove that the functor $M$ constructed above is well behaved with respect to short exact sequences in $\Mod{R}$ and distinguished triangles in $\T$ as well as under the functor $\Ext{}{}(-,-)$.

Note that, in the proof of Theorem \ref{jorgensen_thm}, the assumption that $\Hom{\T}{C}{C}$ is flat as an $R$-module is required for proving the injectivity of a certain map which is used in the construction, see \cite[Lemma 4.5]{Jorgensen2}.

The main result of this paper generalises Theorem \ref{jorgensen_thm} to arbitrary triangulated categories and dispenses with the requirement that $R$ be a principal ideal domain. 
We are also able to prove that the generalised Moore spectra functor $M$ is well behaved with respect to short exact sequences and distinguished triangles as well as under the functor $\Ext{}{}(-,-)$. The next section concerns the construction of generalised Moore spectra.

\section{The general construction}\label{general_construction}

The starting point of the general construction is the following generalisation of a well-known result.

\prop
Let $\T$ be a triangulated category with set indexed coproducts and suppose $C$ is a compact object of $\T$. Let $S=\EndC$. Then the functor $$\Hom{\T}{C}{-}:\Add{C}\rightarrow \Proj{\Sop}$$ is an equivalence of categories.
\label{auslander}
\eprop

\pf
See \cite[Proposition II.2.1]{Auslander} for example. The compactness of $C$ can be used to pass to the infinite additive closure and infinitely generated projective modules.
\epf

\rem
{The fact that $\Hom{\T}{C}{-}:\Add{C}\rightarrow\Proj{\Sop}$ is an equivalence of categories means that it is part of an adjunction which is an equivalence of categories:
$$\xymatrix{\Add{C}\ar@<-1ex>[rr]_-{\Hom{\T}{C}{-}} & & \Proj{\Sop}.\ar@<-1ex>[ll]_-{M_{0}}}$$
In particular, the unit of this adjunction is an isomorphism; see Lemma \ref{unit_isomorphism}.}
\erem

Throughout this paper we shall use the following setup (c.f. Setup \ref{jorgensen_setup}).

\setup
{Let $\T$ be a triangulated category with set indexed coproducts and suppose $C$ is an object of $\T$ satisfying the following assumptions:
\begin{enumerate}
\item $C$ is a compact object of $\T$;
\item Its endomorphism algebra $\Sop=\EndC^{\op}$ has finite global dimension $n$; and,
\item We have $\Hom{\T}{C}{\Sigma^{i}C}=\Hom{\T}{C}{\Sigma^{-i}C}=0$ for $i=1,\ldots,n+1$.
\end{enumerate}}
\label{mainsetup}
\esetup

In \cite{Jorgensen2}, an auxilliary category $\M$, which is a certain full subcategory of $\T$, is introduced; see Definition \ref{jorgensen_M}. We define auxilliary categories $\M_{k}$ for $k\in\mathbb{N}\cup\{0\}$ with a view to arriving at an analogous definition of the auxilliary category $\M$.

\df
{For $0\leq k\leq n-1$, we shall define full subcategories $\M_{k}$ of $\T$ as follows:
$$\M_{k}:=\{X\in\T \ | \ \Hom{\T}{C}{\Sigma^{-i}X}=0 \textnormal{ for } i=1,\ldots, k\}$$
with the convention that $\M_{0}=\T$. Note that $\M_{k}\supseteq\M_{k+1}$.}
\label{myM}
\edf

\thm
Let $\T$ be a triangulated category with set indexed coproducts.  Suppose $C$ is an object of $\T$ satisfying the assumptions of Setup \ref{mainsetup}.
Let $\M=\M_{n}$. Then, the functor
$$\Hom{\T}{C}{-}:\M\rightarrow\Mod{\Sop}$$
has a left adjoint
$$M:\Mod{\Sop}\rightarrow \M.$$
Moreover, the functor $M$ is a full embedding of the module category $\Mod{\Sop}$ into the full subcategory $\M$ of $\T$.
\label{moorespectra}
\ethm

We shall the functor $M:\Mod{\Sop}\rightarrow \M$ obtained in Theorem \ref{moorespectra} as taking values in $\T$ via composition with the inclusion functor $i:\M\hookrightarrow\T$ (c.f. Theorem \ref{jorgensen_thm}).

\cor
Viewing the functor obtained in Theorem \ref{moorespectra} as taking values in $\T$, the functor $M:\Mod{\Sop}\rightarrow \T$ is a full embedding.
\ecor

\section{Proof of Theorem \ref{moorespectra}}\label{proof}

The proof of Theorem \ref{moorespectra} consists of a large induction. For clarity of exposition we have isolated the induction hypotheses below.

\newtheorem{hypos}[dfn]{Hypotheses}
\begin{hypos}
\label{hypotheses}
\textnormal{Under the assumptions of Setup \ref{mainsetup}, for $k\geqslant 0$ we have:
\begin{enumerate}
\item There exists a fully faithful functor $M_{k}:\Projk{k}{\Sop}\rightarrow \T$. \label{hyp1}
\item We have the essential image $\N_{k}=\Essential{M_{k}}$ satisfies $\N_{k}\subseteq\M_{k}$ and \linebreak\mbox{$\Hom{\T}{C}{\N_{k}}\subseteq\Projk{k}{\Sop}$.} \label{hyp2}
\item For $A\in\Projk{k}{\Sop}$ and $X\in\M_{k}$ there is a natural isomorphism
$$\Hom{\Mod{S^{\op}}}{A}{\Hom{\T}{C}{X}}\simeq \Hom{\T}{M_{k}(A)}{X}.$$ (Note that $\Hom{\T}{C}{\M_{k}}$ may not lie in $\Projk{k}{\Sop}$.) \label{hyp3}
\item $\Hom{\T}{C}{\Sigma M_{k}(A)}=\cdots = \Hom{\T}{C}{\Sigma^{n+1-k}M_{k}(A)}=0$ for all $A\in \Projk{k}{\Sop}$. \label{hyp4}
\item $\Hom{\T}{M_{k}(A)}{\Sigma^{-1}X}=0$ for all $A\in \Projk{k}{\Sop}$ and $X\in\M_{k+1}$. \label{hyp5}
\item $M_{k}(P)\in\Add{C}$ for any $P\in \Proj{\Sop}$. \label{hyp6}
\end{enumerate}}
\end{hypos}

\rem
{In condition \eqref{hyp2} of Hypotheses \ref{hypotheses} we take, as a convention, $\N_{0}=\Add{C}$. Note that the fact that $\N_{k}\subseteq \M_{k}$ means we have the following:
\begin{equation}
\Hom{\T}{C}{\Sigma^{-1}M_{k}(A)}=\cdots = \Hom{\T}{C}{\Sigma^{-k}M_{k}(A)}=0
\label{hyp7}
\end{equation} for all $A\in\Projk{k}{\Sop}$.}
\label{hyp7}
\erem

The first three hypotheses contain the content of Theorem \ref{moorespectra} and the last three hypotheses are technical conditions required to prove the first three. In Section \ref{induction_hypotheses} we prove the base step of the induction. In Sections \ref{objects}, \ref{morphisms} and \ref{adjointness} we prove the first three conditions in Hypotheses \ref{hypotheses} before finally verifying the technical conditions in Section \ref{scaffolding}.

\subsection{Base step}\label{induction_hypotheses}

\lemma
Under the assumptions of Setup \ref{mainsetup}, Hypotheses \ref{hypotheses} are true for $k=0$.
\label{basestep}
\elemma

\pf
By convention, we have $\N_{0}=\Add{C}$. From Proposition \ref{auslander}, we have an adjoint equivalence of categories:
$$\xymatrix{\Add{C}=\N_{0}\ar@<-1ex>[rr]_-{\Hom{\T}{C}{-}} & & \Proj{\Sop}.\ar@<-1ex>[ll]_-{M_{0}}}$$
Hence the unit of this adjunction is an isomorphism (see Lemma \ref{unit_isomorphism}), and for $P\in\Proj{\Sop}$ and $X\in \M_{0}=\T$ we have a natural isomorphism
$$\Hom{\Mod{S^{\op}}}{P}{\Hom{\T}{C}{X}}\simeq \Hom{\T}{M_{0}(P)}{X}.$$
(Note that this isomorphism is stronger than just adjointness in this case because $X\in\T$.)
We also have $\N_{0}=\Essential{M_{0}}=\Add{C}\subseteq \M_{0}$. This proves \eqref{hyp1}, \eqref{hyp2} and \eqref{hyp3} of Hypotheses \ref{hypotheses} for $k=0$.

Since $M_{0}(P)\in\Add{C}$ for $P\in\Proj{\Sop}$, the assumption that $\Hom{\T}{C}{\Sigma^{i}C}=0$ for $i=1,\ldots, n+1$ forces
$$\Hom{\T}{C}{\Sigma M_{0}(P)}=\cdots = \Hom{\T}{C}{\Sigma^{n+1}M_{0}(P)}=0.$$
Hence hypothesis \eqref{hyp4} is satisfied for $k=0$.

For $X\in\M_{1}$, we have $\Hom{\T}{M_{0}(P)}{\Sigma^{-1}X}=0$ since $M_{0}(P)\in\Add{C}$, so hypotheses \eqref{hyp5} and \eqref{hyp6} are also satisfied for $k=0$.
This completes the proof of the base step.
\epf

\

We now turn our attention to the proof of the induction step, starting with the construction of the functor $M_{k+1}:\Projk{k+1}{\Sop}\rightarrow\M_{k+1}$. We construct the object map and then describe how $M_{k+1}$ behaves on morphisms.

\subsection{The definition of $M_{k+1}$ on objects}\label{objects}

Let $A\in\Projk{k+1}{\Sop}$ and choose a short exact sequence 
\begin{equation}
\xi^{A}\,:\, 0\rightarrow K^{A}\rightarrow P^{A}\rightarrow A\rightarrow 0.
\label{ses1}
\end{equation}
Then $K^{A}$ and $P^{A}$ are objects in $\Projk{k}{\Sop}$ and, by induction, we have a functor $M_{k}:\Projk{k}{\Sop}\rightarrow\M_{k}$. Applying this functor to the morphism $K^{A}\rightarrow P^{A}$ in $\xi^{A}$ and extending to a distinguished triangle, we obtain
\begin{equation}
\tri{M_{k}(K^{A})}{M_{k}(P^{A})}{M}.
\label{dist_tri1}
\end{equation}
Set $M_{k+1}(A):=M$ in the distinguished triangle \eqref{dist_tri1}.

\lemma
The essential image $\N_{k+1}=\Essential{M_{k+1}}$ satisfies $\N_{k+1}\subseteq\M_{k+1}$. In particular, $M_{k+1}(A)$ is an object of $\M_{k+1}$.
\label{essential_image}
\elemma

\pf
It suffices to show that $M_{k+1}(A)$ as constructed above is an object of $\M_{k+1}$. Applying the functor $\Hom{\T}{C}{-}$ to \eqref{dist_tri1} yields the following long exact sequence:
$$\Hom{\T}{C}{\Sigma^{-i}M_{k}(P^{A})}\rightarrow \Hom{\T}{C}{\Sigma^{-i}M_{k+1}(A)}\rightarrow \Hom{\T}{C}{\Sigma^{-i+1}M_{k}(^{A})}.$$
By Remark \ref{hyp7} (which follows from hypothesis \ref{hypotheses}\eqref{hyp2}), we have $$\Hom{\T}{C}{\Sigma^{-i}M_{k}(K^{A})}=0 \textnormal{ for } i=1,\ldots, k,$$ and by assumption $\Hom{\T}{C}{\Sigma^{-j}C}=0$ for $j=1,\ldots,n+1$. We also have $M_{k}(P)\in\Add{C}$ by hypothesis \eqref{hyp6} of \ref{hypotheses}. Hence, we obtain
$$\Hom{\T}{C}{\Sigma^{-i}M_{k+1}(A)}=0 \textnormal{ for } i=2,\ldots, k+1.$$

By Hypotheses \ref{hypotheses}, we have an adjoint pair
$$\xymatrix{\N_{k}\ar@<-1ex>[rr]_-{\Hom{\T}{C}{-}} & & \Projk{k}{\Sop}\ar@<-1ex>[ll]_-{M_{k}}}$$
whose unit $\eta$ is an isomorphism (since $M_{k}$ is fully faithful). So we have the following commutative square:
$$\xymatrix{\Hom{\T}{C}{M_{k}(K^{A})}\ar[r] & \Hom{\T}{C}{M_{k}(P^{A})} \\
K_{0}^{A}\ar@{^(->}[r]\ar[u]^-{\eta_{K^{A}}}_-{\sim} & P^{A}\ar[u]^-{\eta_{P^{A}}}_-{\sim}}$$
Hence the map $\Hom{\T}{C}{M_{k}(K^{A})}\rightarrow \Hom{\T}{C}{M_{k}(P^{A})}$ is an injection, which forces $\Hom{\T}{C}{\Sigma^{-1}M_{k+1}(A)}=0$. Thus $M_{k+1}(A)\in \M_{k+1}$ and we have $\N_{k+1}\subseteq \M_{k+1}$ and a map $M_{k+1}:\Projk{k+1}{\Sop}\rightarrow \M_{k+1}$.
\epf

\

Lemma \ref{essential_image} means that we now have a map $M_{k+1}:\Projk{k+1}{\Sop}\rightarrow\M_{k+1}$ which will be the object map of the functor $M_{k+1}$.

We next need to explain how $M_{k+1}$ is defined on morphisms.

\subsection{The defintion of $M_{k+1}$ on morphisms}\label{morphisms}

Let $f:A\rightarrow B$ be a morphism in $\Projk{k+1}{\Sop}$ and suppose we have chosen two short exact sequences $\xi^{A}$ and $\xi^{B}$.
Then the morphism $f$ lifts to a commutative diagram of short exact sequences:
$$\xymatrix{\xi^{A}\, : & 0\ar[r] & K^{A}\ar[r]^-{i^{A}}\ar[d]^-{\kappa} & P^{A}\ar[r]^-{p^{A}}\ar[d]^-{\pi} & A\ar[r]\ar[d]^-{f} & 0 \\
			\xi^{B}\, : & 0\ar[r] & K^{B}\ar[r]_-{i^{B}} & P^{B}\ar[r]_-{p^{B}} & B\ar[r] & 0.}$$
By construction of the map $M_{k+1}$, we have the following commutative diagram:
$$\xymatrix{M_{k}(K^{A})\ar[r]^-{M_{k}(i^{A})}\ar[d]^-{M_{k}(\kappa)} & M_{k}(P^{A})\ar[r]\ar[d]^-{M_{k}(\pi)} & M_{k+1}(A)\ar[r]\ar@{-->}[d]^-{\alpha} &\Sigma M_{k}(K^{A})\ar[d]^-{\Sigma M_{k}(\kappa)} \\
M_{k}(K^{B})\ar[r]_-{M_{k}(i^{B})} & M_{k}(P^{B})\ar[r] & M_{k+1}(B)\ar[r] & \Sigma M_{k}(K^{B}).}$$
We claim that $\alpha$ is the unique morphism $M_{k+1}(A)\rightarrow M_{k+1}(B)$ making this diagram commute. Suppose $\alpha':M_{k+1}(A)\rightarrow M_{k+1}(B)$ is another morphism making this diagram commute. Then, by the commutativity of the respective diagrams containing $\alpha$ and $\alpha'$ we obtain a commutative diagram
$$\xymatrix{M_{k}(K^{A})\ar[r] & M_{k}(P^{A})\ar[r]\ar@{-->}[dr]_-{0}\ar[d] & M_{k+1}(A)\ar[r]\ar[d]^{\alpha -\alpha'} & \Sigma M_{k}(K^{A})\ar@{-->}[dl]^{\phi} \\
 & M_{k}(P^{B})\ar[r] & M_{k+1}(B), & }$$
 where the broken arrow $\phi:\Sigma M_{k+1}(K^{A})\rightarrow M_{k+1}(B)$ exists because the composite broken arrow $M_{k}(K^{A})\rightarrow M_{k}(P^{B})$ is zero by the commutative of the diagrams containing the morphisms $\alpha$ and $\alpha'$, respectively. By Hypotheses \ref{hypotheses} \eqref{hyp5}, the broken arrow $\phi=0$, in which case $\alpha=\alpha'$, as claimed.
 
\subsection{Functoriality and adjointness}\label{adjointness}

Sections \ref{objects} and \ref{morphisms} describe how to define $M_{k+1}$ on objects and morphisms, respectively. However, it is still not clear that $M_{k+1}$ is functorial, in particular, that $M_{k+1}(A)$ as described above is independent of the choice of short exact sequence $\xi^{A}$. In this section, we show that $M_{k+1}:\Projk{k+1}{\Sop}\rightarrow \M_{k+1}$ does indeed define a functor and also obtain the adjunction in Hypotheses \ref{hypotheses} via the same classical representability result (stated as Lemma \ref{classical} below).

In the next lemma we first obtain a natural isomorphism on the full subcategory $\M_{k+1}$.

\lemma
For $A\in\Projk{k+1}{\Sop}$ there is a natural isomorphism
$$\Hom{\Mod{\Sop}}{A}{\Hom{\T}{C}{-}}\simeq \Hom{\T}{M}{-}$$
on $\M_{k+1}$.
\label{natural_isomorphism}
\elemma

\pf
Let $A\in\Projk{k+1}{\Sop}$ and choose a short exact sequence $\xi^{A}$ as in \eqref{ses1} and obtain a distinguished triangle as in \eqref{dist_tri1}. 
Let $X,Y\in\M_{k+1}\subseteq\M_{k}$ and apply the functors $\Hom{\Mod{\Sop}}{-}{\Hom{\T}{C}{X}}$ and $\Hom{\T}{-}{X}$ to \eqref{ses1} and \eqref{dist_tri1}, respectively, to obtain the following commutative diagram:
$$\xymatrix{0\ar[r] & (A, (C,X))\ar[r]\ar@{-->}[d]^-{\sim} & (P^{A},(C,X))\ar[r]\ar[d]^-{\sim} & (K^{A},(C,X))\ar[r]\ar[d]^-{\sim} & \\
(M_{k}(K^{A}),\Sigma^{-1}X)\ar[r] & (M, X)\ar[r] & (M_{k}(P^{A}),X)\ar[r] & (M_{k}(K^{A}),X)\ar[r] & }$$
where, in the above diagram, we have used the notation $(-,(C,X))$ and $(-,X)$ as shorthand for the functors $\Hom{\Mod{\Sop}}{-}{\Hom{\T}{C}{X}}$ and $\Hom{\T}{-}{X}$, respectively. By hypothesis \ref{hypotheses}\eqref{hyp5}, $\Hom{\T}{M_{k}(K^{A})}{\Sigma^{-1}X}=0$ for \mbox{$X\in\M_{k+1}$}, hence, the broken arrow exists and is an isomorphism.

$\Hom{\Mod{\Sop}}{-}{\Hom{\T}{C}{X}}$, $\Hom{\Mod{\Sop}}{-}{\Hom{\T}{C}{Y}}$, $\Hom{\T}{-}{X}$ \linebreak and $\Hom{\T}{-}{Y}$ applied to \eqref{ses1} and \eqref{dist_tri1} for $X,Y\in\M_{k+1}$ yield the following commutative diagram:
$$\xymatrix@C-53pt@!C{ &  0\ar[rr] & & (A,(C,Y))\ar[rr]\ar[dd]|\hole & & (P^{A},(C,Y))\ar[dd]|\hole\ar[rr] & &  (K^{A},(C,Y))\ar[dd]|\hole\ar[rr] & & \cdots \\
0\ar[rr] & & (A,(C,X))\ar[ur]\ar[dd]\ar[rr] & & (P^{A},(C,X))\ar[ur]\ar[dd]\ar[rr] & & (K^{A},(C,X))\ar[dd]\ar[ur]\ar[rr] & & \cdots & \\
 & 0\ar[rr]|\hole & & (M, Y)\ar[rr]|\hole & & (M_{k}(P^{A}),Y)\ar[rr]|\hole & & (M_{k}(K^{A}),Y)\ar[rr] & & \cdots \\
0\ar[rr] & & (M,X)\ar[ur]\ar[rr] & & (M_{k}(P^{A}),X)\ar[ur]\ar[rr] & & (M_{k}(K^{A}),X)\ar[ur]\ar[rr] & & \cdots & }$$
and we obtain a natural isomorphism,
$$\Hom{\Mod{\Sop}}{A}{\Hom{\T}{C}{-}}\simeq \Hom{\T}{M}{-},$$
on $\M_{k+1}$, as desired.
\epf

\

We now state the classical representability result which we shall use.

\lemma
[\cite{MacLane}, Corollary IV.1.2]
Let $\C$ and $\D$ be categories. A functor $G:\C\rightarrow \D$ has a left adjoint if and only if for every object $C$ of $\C$ there is a natural isomorphism
$$\phi:\Hom{\C}{C}{G(D)}\simeq \Hom{\D}{F_{0}(C)}{D}$$
which is natural in $D\in\D$. Then $F_{0}$ is the object function of the left adjoint of $G$.
\label{classical}
\elemma

Given a functor $M_{k}:\Projk{k}{\Sop}\rightarrow \M_{k}$ we have constructed a map $M_{k+1}:\Projk{k+1}{\Sop}\rightarrow\M_{k+1}$. Via Lemma \ref{classical}, the natural isomorphism obtained in Lemma \ref{natural_isomorphism} would mean that $M_{k+1}$ is left adjoint to $\Hom{\T}{C}{-}$, (and, in particular, is a functor) provided that $\Hom{\T}{C}{\M_{k+1}}$ lies in $\Projk{k+1}{\Sop}$. However, this may not be the case and the rest of this section is devoted to establishing this fact.

It is sufficient to prove that the functor $\Hom{\T}{C}{-}$ when it acts on the essential image $\N_{k+1}$ (which is contained in $\M_{k+1}$ by Lemma \ref{essential_image}) takes values in $\Projk{k+1}{\Sop}$, then we obtain an adjunction:
$$\xymatrix{\N_{k+1}\ar@<-1ex>[rr]_-{\Hom{\T}{C}{-}} & & \Projk{k+1}{\Sop}.\ar@<-1ex>[ll]_-{M_{k+1}}}$$
Composition of the functor $M_{k+1}:\Projk{k+1}{\Sop}\rightarrow \N_{k+1}$ with the inclusion functor $\iota: \N_{k+1}\rightarrow \M_{k+1}$ gives the desired functor. Thus we need to prove the following lemma.

\lemma
The restriction of the functor $\Hom{\T}{C}{-}:\T\rightarrow \Mod{\Sop}$ to the full subcategory $\N_{k+1}$ of $\T$ takes values in $\Projk{k+1}{\Sop}$ and hence defines a functor
$$\Hom{\T}{C}{-}:\N_{k+1}\rightarrow\Projk{k+1}{\Sop}.$$
\label{restriction_functor}
\elemma

\pf
Suppose $X\in\N_{k+1}$, then there is an $\Sop$-module $A$ of projective dimension at most $k+1$ such that $M_{k+1}(A)\cong X$. 
Take the usual short exact sequence $\xi^{A}$ (see \eqref{ses1}) and obtain the usual distinguished triangle \eqref{dist_tri1}.
Applying the functor $\Hom{\T}{C}{-}$ to \eqref{dist_tri1} gives the following long exact sequence:
$$0\rightarrow (C, M_{k}(K^{A}))\rightarrow (C,M_{k}(P^{A}))\rightarrow (C, M_{k+1}(A))\rightarrow (C, \Sigma M_{k}(K^{A})),$$
where we have used the shorthand described in the proof of Lemma \ref{natural_isomorphism} to denote the Hom-spaces. The zero on the left hand side comes by the fact that $\Hom{\T}{C}{\Sigma^{-1}M_{k+1}(A)}=0$ because $M_{k+1}(A)\in\M_{k+1}$ (hypothesis \eqref{hyp5} of \ref{hypotheses}). By condition \eqref{hyp4} of Hypotheses \ref{hypotheses}, we have $\Hom{\T}{C}{\Sigma M_{k}(K^{A})}=0$. Hence, we have the commutative diagram below.
$$\xymatrix{0\ar[r] & (C, M_{k}(K^{A}))\ar[r] & (C,M_{k}(P^{A}))\ar[r] & (C, M_{k+1}(A))\ar[r] & 0. \\
 & K^{A}\ar[u]^-{\eta_{K^{A}}}_-{\sim}\ar@{^(->}[r] & P^{A}\ar[u]^-{\eta_{P^{A}}}_-{\sim} & & }$$
Therefore, $\Hom{\T}{C}{M_{k+1}(A)}\in\Projk{k+1}{\Sop}$. Moreover, $\Hom{\T}{C}{M_{k+1}(A)}\cong A$ for all $A\in\Projk{k+1}{\Sop}$. It follows that, $\Hom{\T}{C}{X}\in\Projk{k+1}{\Sop}$.
\epf

\cor
The construction of sections \ref{objects} and \ref{morphisms} defines a functor $M_{k+1}:\Projk{k+1}{\Sop}\rightarrow\M_{k+1}$. Moreover, $M_{k+1}$ is fully faithful.
\ecor

\pf
Observe that in the proof of Lemma \ref{restriction_functor} we obtained the isomorphism $\Hom{\T}{C}{M_{k+1}(A)}\cong A$. It follows by Lemma \ref{unit_isomorphism} that $M_{k+1}$ is fully faithful.
\epf

\

We now summarise the work of sections \ref{objects}, \ref{morphisms} and \ref{adjointness} with the following proposition.

\prop
Under the assumptions of Setup \ref{mainsetup} and the hypotheses of Hypotheses \ref{hypotheses}, we have:
\begin{enumerate}
\item There exists a fully faithful functor $M_{k+1}:\Projk{k+1}{\Sop}\rightarrow \M_{k+1}$.
\item The essential image $\N_{k+1}=\Essential{M_{k+1}}$ satisfies $\N_{k+1}\subseteq\M_{k+1}$ and 
\linebreak\mbox{$\Hom{\T}{C}{\N_{k+1}}\subseteq\Projk{k+1}{\Sop}.$}
\item For $A\in\Projk{k+1}{\Sop}$ and $X\in\M_{k+1}$ there is a natural isomorphism
$$\Hom{\Mod{S^{\op}}}{A}{\Hom{\T}{C}{X}}\simeq \Hom{\T}{M_{k+1}(A)}{X}.$$
\end{enumerate}
\label{firstthree}
\eprop

\rem
{In Proposition \ref{firstthree} we have accomplished the proof of the induction step for the first three hypotheses in Hypotheses \ref{hypotheses}. However, in obtaining this proof we have used all the hypotheses stated in Hypotheses \ref{hypotheses}. The proofs of these technical hypotheses are straightforward and are explained in the next section.}
\erem

\subsection{The scaffolding}\label{scaffolding}

We start the construction of the scaffolding with the easier hypotheses, \eqref{hyp4}, \eqref{hyp5} and \eqref{hyp6} of Hypotheses \ref{hypotheses}, dealing with each in turn.

\lemma
Under the assumptions of Setup \ref{mainsetup} and Hypotheses \ref{hypotheses} we have:
$$\Hom{\T}{C}{\Sigma M_{k+1}(A)}=\cdots = \Hom{\T}{C}{\Sigma^{n-k}M_{k+1}(A)}=0$$
for all $A\in\Projk{k+1}{\Sop}$.
\label{lemma_hyp4}
\elemma

\pf
Suppose $A\in\Projk{k+1}{\Sop}$ and consider the distinguished triangle \eqref{dist_tri1} obtained from the short exact sequence $\xi^{A}$ (see \eqref{ses1}). Note that in $\xi^{A}$ the $\Sop$-module $K^{A}$ is an object of $\Projk{k}{\Sop}$ and $P^{A}$ is projective. Apply the functor $\Hom{\T}{C}{-}$ to \eqref{dist_tri1}, then the claim of the lemma can be read off from the resulting long exact sequence after making the observation that by Hypotheses \ref{hypotheses}, $\Hom{\T}{C}{\Sigma^{i}M_{k}(K^{A})}=0$ for $i=1,\ldots, n+1-k$, and $M_{k}(P^{A})\in\Add{C}$ so that $\Hom{\T}{C}{\Sigma^{i}M_{k}(P^{A})}=0$ for $i=1,\ldots, n+1$. Condition \eqref{hyp4} of Hypotheses \ref{hypotheses} now follows for $k+1$.
\epf

\lemma
Under the assumptions of Setup \ref{mainsetup} and Hypotheses \ref{hypotheses} we have $\Hom{\T}{M_{k+1}(A)}{\Sigma^{-1}X}=0$ for all $A\in\Projk{k+1}{\Sop}$ and $X\in\M_{k+2}$.
\label{lemma_hyp5}
\elemma

\pf
Let $A\in\Projk{k+1}{\Sop}$ and $X\in\M_{k+2}$. Consider the distinguished triangle \eqref{dist_tri1} obtained from the short exact sequence $\xi^{A}$. 
Applying the functor $\Hom{\T}{-}{X}$ to distinguished triangle \eqref{dist_tri1} with $X\in\M_{k+2}$. We get:
$$\Hom{\T}{\Sigma^{2}M_{k}(K^{A})}{X}\rightarrow \Hom{\T}{\Sigma M_{k+1}(A)}{X}\rightarrow \Hom{\T}{\Sigma M_{k}(P^{A})}{X}.$$

We know that $\Hom{\T}{\Sigma M_{k}(P^{A})}{X}=0$ because $X\in\M_{k+2}$ and, by hypothesis \eqref{hyp6} of \ref{hypotheses}, we have $M_{k}(P^{A})\in\Add{C}$. So, we only need to see that $\Hom{\T}{\Sigma^{2}M_{k}(K^{A})}{X}=0$. 

Choose short exact sequences
$$\xi_{i}^{A}\,:\, 0\rightarrow K_{i}^{{A}}\rightarrow P_{i}^{A}\rightarrow K_{i-1}^{A}\rightarrow 0,$$
for $i\geqslant 1$, where we take the convention $K_{0}^{A}=K^{A}$ and $\xi_{0}^{A}=\xi^{A}$. (For instance, take each short exact sequence $\xi_{i}^{A}$ to be the $i^{\textnormal{th}}$-stage of a projective resolution of $A$ of shortest length.)
By the same argument, $\Hom{\T}{\Sigma^{2}M_{k}(K^{A})}{X}=0$ if $\Hom{\T}{\Sigma^{3}M_{k}(K_{1}^{A})}{X}=0$. Applying the same argument inductively yields the following implication:
$$\Hom{\T}{\Sigma^{i+3}M_{k}(K_{i+1}^{A})}{X}=0 \implies \Hom{\T}{\Sigma^{i+2}M_{k}(K_{i}^{A})}{X}=0.$$
Eventually, $K^{A}_{k}=P^{A}_{k+1}$ for some projective $P^{A}_{k+1}$ and $X\in\M_{k+2}$, so
$$\Hom{\T}{\Sigma^{k+2}M_{k}(K_{k}^{A})}{X}=\Hom{\T}{\Sigma^{k+2}M_{k}(P_{k+1}^{A})}{X}=0,$$
because $M_{k}(P_{k+1}^{A})\in\Add{C}$,
so it follows that $\Hom{\T}{\Sigma^{k+1}M_{k}(K_{k-1}^{A})}{X}=0$, and eventually we see that $$\Hom{\T}{\Sigma M_{k+1}(A)}{X}=0.$$ That is, for $A\in\Projk{k+1}{\Sop}$ and $X\in\M_{k+2}$, we have $\Hom{\T}{M_{k+1}(A)}{\Sigma^{-1}X}=0$, as required.
\epf

\lemma
Under the assumptions of Setup \ref{mainsetup} and Hypotheses \ref{hypotheses} then given any projective $\Sop$-module $P$ we have  $M_{k+1}(P)\in \Add{C}$.
\label{lemma_hyp6}
\elemma

\pf
Let $P$ be a projective $\Sop$-module and consider the short exact sequence coming from its projective resolution:
$$0\longrightarrow 0 \longrightarrow P\rightiso P\longrightarrow 0.$$
This gives a distinguished triangle:
$$0\rightarrow M_{k}(P)\rightarrow M_{k+1}(P)\rightarrow \Sigma 0,$$
by construction. Hence $M_{k+1}(P)\cong M_{k}(P)$ and $M_{k}(P)\in\Add{C}$ by induction, so we have $M_{k+1}(P)\in\Add{C}$.
\epf

\

Hence we have shown that if Hypotheses \ref{hypotheses} are true for $k$ then they are true for $k+1$. We know, by Lemma \ref{basestep} that they are true for $k=0$.
It is clear that the induction terminates at $k=n$, thus setting $M=M_{n}$ and $\M=\M_{n}$ gives the existence of a functor
$$M:\Mod{\Sop}\rightarrow \M$$
which is left adjoint to the functor $\Hom{\T}{C}{-}:\M\rightarrow \Mod{\Sop}$, completing the proof of Theorem \ref{moorespectra}.

\

It is useful for the next sections to highlight an important aspect of the proof of Theorem \ref{moorespectra}. 

\rem
{Given a short exact sequence, $0\rightarrow K^{A}\rightarrow P^{A}\rightarrow A\rightarrow 0$, in $\Mod{\Sop}$ with $K^{A}\in\Projk{k}{\Sop}$, $P^{A}\in\Proj{\Sop}$ and $A\in\Projk{k+1}{\Sop}$, the proof of Theorem \ref{moorespectra} gives a distinguished triangle $\tri{MK}{MP}{MA}$.}
\label{remark0}
\erem

\section{J\o rgensen's construction revisited}\label{revisited}

In this brief section we return to J\o rgensen's theorem (stated in this paper as Theorem \ref{jorgensen_thm}). We first note the following specialisation of Theorem \ref{moorespectra}.

\prop
Let $\T$ be an $R$-linear triangulated category with set indexed coproducts.  Suppose $C$ is a compact object of $\T$ satisfying the following assumptions:
\begin{enumerate}
\item Its endomorphism algebra $\Sop=\EndC^{\op}$ has global dimenison $1$; and,
\item We have $\Hom{\T}{C}{\Sigma C}=\Hom{\T}{C}{\Sigma^{-1}C}=0$.
\end{enumerate}
Let $\M=\M_{1}$ (see Definition \ref{myM}). Then, the functor
$$\Hom{\T}{C}{-}:\M\rightarrow\Mod{\Sop}$$
has a left adjoint
$$M:\Mod{\Sop}\rightarrow \M.$$
Moreover, the functor $M$ is a full embedding of the module category $\Mod{\Sop}$ into the full subcategory $\M$ of $\T$.
\label{moorespectra_special}
\eprop

Proposition \ref{moorespectra_special}  can be viewed as an analogue of Theorem \ref{jorgensen_thm}, and hence Theorem \ref{moorespectra} can be considered a higher analogue of Theorem \ref{jorgensen_thm}. Indeed, the hypotheses of Setup \ref{jorgensen_setup} imply the hypotheses of Proposition \ref{moorespectra_special}. In addition, if the unit of the adjunction obtained in Theorem \ref{jorgensen_thm} is an isomorphism then we have
$$R\cong\Hom{\T}{C}{M'(R)}\cong\Hom{\T}{C}{C}=S,$$
where $M'$ denotes the functor obtained in Theorem \ref{jorgensen_thm}. Hence, in this case, $\Mod{R^{\op}}\simeq\Mod{\Sop}$ and the functors $M$, obtained in Proposition \ref{moorespectra_special}, and $M'$, obtained in Theorem \ref{jorgensen_thm}, coincide.

\rem
{Note that the full subcategory $\N_{1}$, defined in Hypotheses \ref{hypotheses}, of the auxilliary category $\M_{1}$ defined in Definition \ref{myM} coincides with the auxilliary category $\M$ defined in Definition \ref{jorgensen_M}.}
\erem 

\section{The functor $M$ is well behaved}\label{well_behaved}

We now show that the functor $M:\Mod{\Sop}\rightarrow\T$ constructed in Theorem \ref{moorespectra} is well behaved with respect to short exact sequences in $\Mod{\Sop}$ and distinguished triangles in $\T$. The hard work carried out in section \ref{general_construction} provides the setting for a functorial proof of Theorem \ref{sequence_triangle} which differs in character entirely with the corresponding result which it generalises (\cite[Theorem 4.11]{Jorgensen2}).

\thm
Let $\T$ be a triangulated category with set indexed coproducts. Let $C$ be an  object of $\T$ satisfying the assumptions of Setup \ref{mainsetup}. Let $M:\Mod{\Sop}\rightarrow \T$ be the functor obtained in Theorem \ref{moorespectra}. If $0\longrightarrow A'\rightlabel{a'} A\rightlabel{a}A''\longrightarrow 0$ is a short exact sequence in $\Mod{\Sop}$, then there is a distinguished triangle:
$$\trilabels{MA'}{MA}{MA''}{Ma'}{Ma}{}$$
in $\T$.
\label{sequence_triangle}
\ethm

\pf
We prove the result for each functor $M_{k}:\Projk{k}{\Sop}\rightarrow \M_{k}$ by induction on $k$. Suppose $k=0$ and suppose we have a short exact sequence $0\longrightarrow P'\rightlabel{\pi'}P\rightlabel{\pi}P''\longrightarrow 0$ of projective $\Sop$-modules. Such an exact sequence is split, hence there are splitting maps:
$$\xymatrix{0\ar[r] & P'\ar@<1ex>[r]^-{\pi'} & P\ar@<1ex>[r]^-{\pi}\ar@<1ex>[l]^-{p} & P''\ar[r]\ar@<1ex>[l]^{p''} & 0}$$
such that $\textnormal{id}_{P}=p''\pi+\pi'p$. Applying the functor $M_{0}$ to the diagram above yields:
$$\xymatrix{M_{0}(P')\ar@<1ex>[r]^-{M_{0}(\pi')} & M_{0}(P)\ar@<1ex>[r]^-{M_{0}(\pi)}\ar@<1ex>[l]^-{M_{0}(p)} & M_{0}(P'').\ar@<1ex>[l]^{M_{0}(p'')}}$$
It is well known that such a diagram is isomorphic to a distinguished triangle:
$$\tri{M_{0}(P')}{M_{0}(P')\coprod M_{0}(P'')}{M_{0}(P'')}.$$
Hence, $\trilabels{M_{0}(P')}{M_{0}(P)}{M_{0}(P'')}{M_{0}(\pi')}{M_{0}(\pi)}{}$ is a distinguished triangle, proving the assertion for $k=0$.

Let $k\geq 1$ and, by induction, suppose that any short exact sequence $0\longrightarrow K'\rightlabel{\kappa'} K\rightlabel{\kappa} K''\longrightarrow 0$ in $\Projk{k}{\Sop}$ corresponds to a distinguished triangle $\trilabels{M_{k}(K')}{M_{k}(K)}{M_{k}(K'')}{M_{k}(\kappa')}{M_{k}(\kappa)}{}$. Now suppose we have a short exact sequence $0\longrightarrow A'\rightlabel{a'} A\rightlabel{b} A''\longrightarrow 0$ with $A',A,A''\in\Projk{k+1}{\Sop}$. By \cite[Horseshoe Lemma 2.2.8]{Weibel}, one can obtain a diagram of the form:
$$\xymatrix{ & & 0\ar[d] & 0\ar[d] & 0\ar[d] & \\
\xi^{A'}\,: & 0\ar[r] & K'\ar[r]^-{f'}\ar[d]^-{\kappa'} & P'\ar[r]^-{g'}\ar[d]^-{\pi'} & A'\ar[r]\ar[d]^-{a'} & 0 \\
\xi^{A}\,: & 0\ar[r]	& K\ar[r]^-{f}\ar[d]^-{\kappa}  & P\ar[r]^-{g}\ar[d]^-{\pi}	  & A\ar[d]^-{a}\ar[r] &  0\\
\xi^{A''}\,: & 0\ar[r] & K''\ar[r]^-{f''}\ar[d] & P''\ar[r]^-{g''}\ar[d] & A''\ar[r]\ar[d] & 0 \\
	&	& 	0	   &		0	& 0 & }$$
where $P=P'\coprod P''$ and $K\in\Projk{k}{\Sop}$. By induction, the short exact sequence in the left hand column gives $\trilabels{M_{k}(K')}{M_{k}(K)}{M_{k}(K'')}{M_{k}(\kappa')}{M_{k}(\kappa)}{}$ and the short exact sequence of projectives sitting in the central column gives a split distinguished triangle $\trilabels{M_{k}(P')}{M_{k}(P)}{M_{k}(P'')}{M_{k}(\pi')}{M_{k}(\pi)}{}$. By construction, we obtain the following diagram:
\begin{equation*}
\xymatrix{M_{k}(K')\ar[r]^-{M_{k}(f')}\ar[d]_-{M_{k}(\kappa')} & M_{k}(P')\ar[r]^-{M_{k}(g')}\ar[d]^-{M_{k}(\pi')} & M_{k+1}(A')\ar[r]^-{h'}\ar[d]^-{M_{k+1}(a')} & \Sigma M_{k}(K')\ar[d]^-{\Sigma M_{k}(\kappa')} \\
M_{k}(K)\ar[r]^-{M_{k}(f)}\ar[d]_-{M_{k}(\kappa)} & M_{k}(P)\ar[r]^-{M_{k}(g)}\ar[d]^-{M_{k}(\pi)} & M_{k+1}(A)\ar[r]^-{h}\ar[d]^-{M_{k+1}(a)} & \Sigma M_{k}(K)\ar[d]^-{\Sigma M_{k}(\kappa)} \\
M_{k}(K'')\ar[r]^-{M_{k}(f'')}\ar[d] & M_{k}(P'')\ar[r]^-{M_{k}(g'')}\ar[d] & M_{k+1}(A'')\ar[r]^-{h''}\ar@{-->}[d] & \Sigma M_{k}(K'')\ar[d] \\
\Sigma M_{k}(K')\ar[r]^-{\Sigma M_{k}(f')} & \Sigma M_{k}(P')\ar[r]^-{\Sigma M_{k}(g')} & \Sigma M_{k+1}(A')\ar[r]^-{-\Sigma h'} & \Sigma^{2}M_{k}(K'),}
\label{bigsquare1}
\end{equation*}
which is commutative except for the bottom right hand square, which is anti-commutative. The rows are distinguished triangles (Remark \ref{remark0}), and the broken arrow $MB''\rightarrow \Sigma MB'$ exists by virtue of the axioms of triangulated categories. Now $M_{k+1}(A^{?})$ is constructed as the mapping cone of the map $M_{k}(K^{?})\rightlabel{M_{k}(f^{?})} M_{k}(P^{?})$, where $?$ is either $'$, empty, or $''$; see Section \ref{objects}. By a $3\times 3$ lemma for triangulated categories (see \cite[Lemma 1.7]{Neeman_axioms}, for instance), it follows that the third column is a distinguished triangle. In particular, when $k=n$, one obtains that $\trilabels{MA'}{MA}{MA''}{Ma'}{Ma}{}$ is a distinguished triangle.
\epf

\rem
{Note that the third morphism in the triangle $\trilabels{MA'}{MA}{MA''}{Ma'}{Ma}{}$ in Theorem \ref{sequence_triangle} is unique since, by the proof of Theorem \ref{moorespectra} in Section \ref{proof}, we have $\Hom{\T}{\Sigma MA'}{MA''}=0$.}
\erem

We now aim to prove that for $A,B\in\Mod{\Sop}$ there are natural maps $$\xymatrix{\Ext{n}{\Sop}(A,B)\ar[rr]^-{\Delta^{n}_{A,B}} & & \Hom{\T}{MA}{\Sigma^{n}MB}.}$$
We need to appeal to the definitions of a $\delta$-functor, a universal $\delta$-functor and the fact that the functor $\Ext{}{}(-,-)$ is a universal $\delta$-functor. The following definitions are taken from \cite{Hartshorne1}.

\df
{Let $\A$ and $\B$ be abelian categories. A \textit{(covariant) $\delta$-functor} from $\A$ to $\B$ is a collection of functors $T=(T^{i})_{i\geqslant 0}$ together with a morphism $$\delta^{i}:T^{i}(A'')\rightarrow T^{i+1}(A')$$ for each short exact sequence $0\rightarrow A'\rightarrow A\rightarrow A''\rightarrow 0$ and each $i\geqslant 0$, such that
\begin{enumerate}
\item For each short exact sequence, as above, there is a long exact sequence
$$0\longrightarrow T^{0}(A')\longrightarrow T^{0}(A)\longrightarrow T^{0}(A'')\rightlabel{\delta^{0}} T^{1}(A')\longrightarrow\cdots$$
$$\cdots\longrightarrow T^{i}(A')\longrightarrow T^{i}(A)\longrightarrow T^{i}(A'')\rightlabel{\delta^{i}} T^{i+1}(A')\longrightarrow\cdots$$
\item For each morphism of one short exact sequence, as above, into another $0\rightarrow B'\rightarrow B\rightarrow B''\rightarrow 0$, the $\delta$s give a commutative diagram
$$\xymatrix{T^{i}(A'')\ar[r]^-{\delta^{i}}\ar[d] & T^{i+1}(A')\ar[d] \\
T^{i}(B'')\ar[r]^{\delta^{i}} & T^{i+1}(B').}$$
\end{enumerate}
A \textit{contravariant $\delta$-functor} is defined similarly.}
\edf

\df
{A $\delta$-functor $T=(T^{i})_{i\geqslant 0}:\A\rightarrow\B$ is said to be a \textit{universal $\delta$-functor} if, given any other $\delta$-functor $U=(U^{i})_{i\geqslant 0}:\A\rightarrow\B$ and any given morphism of functors $f^{0}: T^{0}\rightarrow U^{0}$, there exists a unique sequence of morphisms $f^{i}:T^{i}\rightarrow U^{i}$ for each $i\geqslant 0$, starting with the given $f$, which commute with the $\delta^{i}$s for each short exact sequence.}
\edf

It is a well known fact that $\Ext{n}{}(A,-)$ is a covariant universal $\delta$-functor and $\Ext{n}{}(-,B)$ is a contravariant universal $\delta$-functor. We will need the following lemma.

\lemma
Let $\T$ be a triangulated category with set indexed coproducts. Let $C$ be an object of $\T$ satisfying the assumptions of Setup \ref{mainsetup}. Recall that $\Sop=\Endo{C}^{\op}$. Then, we have the following:
\begin{enumerateroman}
\item For $A\in\Mod{\Sop}$ the functor $(\Hom{\T}{MA}{\Sigma^{n}M(-)})_{n\geqslant 0}$ is a covariant $\delta$-functor.
\item For $B\in\Mod{\Sop}$ the functor $(\Hom{\T}{M(-)}{\Sigma^{n}MB})_{n\geqslant 0}$ is a contravariant $\delta$-functor.
\end{enumerateroman}
\label{moore_lemma1}
\elemma

\pf
Let $A\in\Mod{\Sop}$ and consider $U^{n}(-)=\Hom{\T}{MA}{\Sigma^{n}M(-)}$. We claim that $(U^{n})_{n\geqslant 0}$ is a covariant $\delta$-functor. Suppose we have a short exact sequence $0\rightarrow B'\rightarrow B\rightarrow B''\rightarrow 0$ in $\Mod{\Sop}$. By Theorem \ref{sequence_triangle}, there is a distinguished triangle $\tri{MB'}{MB}{MB''}$. This distinguished triangle induces a long exact sequence
$$\cdots\rightarrow U^{i}(B')\rightarrow U^{i}(B)\rightarrow U^{i}(B'')\rightarrow U^{i+1}(B')\rightarrow\cdots$$
and gives a morphism $\delta^{i}_{B}:U^{i}(B'')\rightarrow U^{i+1}(B')$.

Now suppose we have another short exact sequence $0\rightarrow C'\rightarrow C\rightarrow C''\rightarrow 0$ and a morphism of short exact sequences:
$$\xymatrix{0\ar[r] & B'\ar[r]\ar[d] & B\ar[r]\ar[d] & B''\ar[r]\ar[d] & 0 \\
0\ar[r] & C'\ar[r] & C\ar[r] & C''\ar[r] & 0,}$$
which gives a commutative diagram of distinguished triangles:
$$\xymatrix{MB'\ar[r]\ar[d] & MB\ar[r]\ar[d] & MB''\ar[r]\ar[d] & \Sigma MB'\ar[d] \\
MC'\ar[r] & MC\ar[r] & MC''\ar[r] & \Sigma MC'.}$$
This, in turn, yields a commutative diagram:
$$\xymatrix{U^{i}(B'')\ar[r]^{\delta^{i}_{B}}\ar[d] & U^{i+1}(B')\ar[d] \\
U^{i}(C'')\ar[r]_-{\delta^{i}_{C}} & U^{i+1}(C').}$$
Hence, $(U^{n})_{n\geqslant 0}$ is a covariant $\delta$-functor. This proves assertion (i); assertion (ii) is proved similarly.
\epf

\thm
Let $\T$ be a triangulated category with set indexed coproducts. Let $C$ be an object of $\T$ satisfying the assumptions of Setup \ref{mainsetup}. Recall that $\Sop=\Endo{C}^{\op}$. Then, we have the following:
\begin{enumerateroman}
\item For $A,B\in\Mod{\Sop}$ there exist maps, $$\xymatrix{\Ext{n}{\Sop}(A,B)\ar[rr]^-{\Delta^{n}_{A,B}} & & \Hom{\T}{MA}{\Sigma^{n}MB},}$$ which are natural in $A$ and $B$.
\item Given a short exact sequence $0\rightarrow A'\rightarrow A\rightarrow A''\rightarrow 0$ in $\Mod{\Sop}$, we have the following commutative diagram:
$$\xymatrix{0\ar[r] &\Hom{\Sop}{A''}{B}\ar[r]\ar[d]^-{\Delta^{0}_{A'',B}} & \Hom{\Sop}{A}{B}\ar[r]\ar[d]^-{\Delta^{0}_{A,B}} & \Hom{\Sop}{A'}{B}\ar[r]\ar[d]^-{\Delta^{0}_{A',B}} & \\
0\ar[r] & (MA'',MB)\ar[r] & (MA,MB)\ar[r] & (MA',MB)\ar[r] & }$$
$$\xymatrix{\cdots\ar[r] &\Ext{n}{\Sop}(A'',B)\ar[r]\ar[d]^-{\Delta^{n}_{A'',B}} & \Ext{n}{\Sop}(A,B)\ar[r]\ar[d]^-{\Delta^{n}_{A,B}} & \Ext{n}{\Sop}(A',B)\ar[r]\ar[d]^-{\Delta^{n}_{A',B}} & \cdots \\
\cdots\ar[r] & (MA'',\Sigma^{n}MB)\ar[r] & (MA,\Sigma^{n}MB)\ar[r] & (MA',\Sigma^{n}MB)\ar[r] & \cdots. }$$
\item Given a short exact sequence $0\rightarrow B'\rightarrow B\rightarrow B''\rightarrow 0$ in $\Mod{\Sop}$, we have the following commutative diagram:
$$\xymatrix{0\ar[r] &\Hom{\Sop}{A}{B'}\ar[r]\ar[d]^-{\Delta^{0}_{A,B'}} & \Hom{\Sop}{A}{B}\ar[r]\ar[d]^-{\Delta^{0}_{A,B}} & \Hom{\Sop}{A}{B''}\ar[r]\ar[d]^-{\Delta^{0}_{A,B''}} & \\
0\ar[r] & (MA,MB')\ar[r] & (MA,MB)\ar[r] & (MA,MB'')\ar[r] & }$$
$$\xymatrix{\cdots\ar[r] &\Ext{n}{\Sop}(A,B')\ar[r]\ar[d]^-{\Delta^{n}_{A,B}} & \Ext{n}{\Sop}(A,B)\ar[r]\ar[d]^-{\Delta^{n}_{A,B}} & \Ext{n}{\Sop}(A,B'')\ar[r]\ar[d]^-{\Delta^{n}_{A,B''}} & \cdots \\
\cdots\ar[r] & (MA,\Sigma^{n}MB')\ar[r] & (MA,\Sigma^{n}MB)\ar[r] & (MA,\Sigma^{n}MB'')\ar[r] & \cdots.}$$
\end{enumerateroman}
Note that the use of $(X,\Sigma^{n}Y)$ on the bottom row in statements (ii) and (iii) is shorthand for $\Hom{\T}{X}{\Sigma^{n}Y}$.
\label{natural_maps}
\ethm

\pf
The functors $(\Ext{n}{\Sop}(A,-))_{n\geqslant 0}$ and $(\Ext{n}{\Sop}(-,B))_{n\geqslant 0}$ are universal $\delta$-functors. We also know that the functor $(\Hom{\T}{MA}{\Sigma^{n}M(-)})_{n\geqslant 0}$ is a covariant $\delta$-functor and the functor $(\Hom{\T}{-}{\Sigma^{n}MB})_{n\geqslant 0}$ is a contravariant $\delta$-functor by Lemma \ref{moore_lemma1}; the theorem now follows.
\epf

\section{An example from $u$-cluster categories}\label{cluster_categories}

In this section we shall consider a special case of Theorem \ref{moorespectra} when the endomorphism ring of the compact object $C$ of the triangulated category $\T$ is right coherent. In particular, this allows us to specialise Theorem \ref{moorespectra} to the category $\mod{\Sop}$, the full subcategory of $\Mod{\Sop}$ consisting of finitely presented $\Sop$-modules. We shall then apply this specialisation to the case of a path algebra of a quiver which has no oriented cycles. Such a path algebra is well known to be hereditary, and as such coherent, see \cite{Lam} for example. When applying this specialisation of the main theorem to this case, the full embedding of the theorem recovers the canonical embedding of the module category into its $u$-cluster category, where $u\geqslant 2$ is an integer.

\subsection{A version of the main result for finitely presented modules}

We first specialise the main result to the case for finitely presented $\Sop$-modules. For the basic facts on coherent rings and modules we refer to \cite{Lam}.

\df
{A ring $R$ is said to be \textit{right coherent} if every finitely generated right ideal of $R$ is also finitely presented as a right $R$-module.}

\textnormal{A finitely generated right $R$-module $A$ is said to be \textit{coherent} if every finitely generated submodule of $A$ is finitely presented.}
\edf

Recall that a ring $R$ is right coherent if and only if any finitely presented right $R$-module is coherent. It follows that, if $R$ is right coherent then every finitely generated projective right $R$-module is coherent. It is well known that the kernel of a homomorphism of finitely generated projective $R$-modules is also finitely generated, see \cite[Lemma 2.11]{Holm_Jorgensen1}. Hence, the kernel of a homomorphism of finitely generated right $R$-modules is finitely presented. The following is now an easy lemma.

\lemma
Suppose $R$ is a right coherent ring. A finitely presented right $R$-module $A$ with finite projective dimension $k$ has a projective resolution of length $k$ consisting of finitely generated projective right $R$-modules.
\label{coherent5}
\elemma

In light of Lemma \ref{coherent5} and the usual finite version of Proposition \ref{auslander}, we now obtain the following version of Theorem \ref{moorespectra}.

\thm
Let $\T$ be a triangulated category with set indexed coproducts.  Suppose $C$ is an object of $\T$ satisfying the following assumptions:
\begin{enumerate}
\item Its endomorphism algebra $\Sop=\EndC^{\op}$ has finite global dimension $n$;
\item We have $\Hom{\T}{C}{\Sigma^{i}C}=\Hom{\T}{C}{\Sigma^{-i}C}=0$ for $i=1,\ldots, n+1$; and,
\item The endomorphism algebra $S$ is right coherent. 
\end{enumerate}
Let $\M=\M_{n}$, then there exists a full embedding
$M:\mod{\Sop}\rightarrow \M$.
\label{moorespectrafinite}
\ethm

Recall that the category $\mod{\Sop}$ of finitely presented right $S$-modules is an abelian category if and only if $S$ is right coherent. Thus, given Lemma \ref{coherent5}, it follows that the proofs of Theorems \ref{sequence_triangle} and \ref{natural_maps} can be used to prove versions of these theorems for finitely presented right $S$-modules.

\rem
{Note that in Theorem \ref{moorespectrafinite} we do not obtain that $M$ is left adjoint to $\Hom{\T}{C}{-}$ because it is not clear that $\Hom{\T}{C}{-}$ applied to $\M$ necessarily takes values in $\mod{\Sop}$. The construction of Theorem \ref{moorespectra} applies to Theorem \ref{moorespectrafinite} by virtue of Lemma \ref{coherent5} and the fact that one can still show at each stage of the construction that there is an adjunction
$$\xymatrix{\N_{k}\ar@<-1ex>[rr]_-{\Hom{\T}{C}{-}} & & \projk{k}{\Sop},\ar@<-1ex>[ll]_-{M_{k}}}$$
where $\N_{k}=\{X\in \T \ | \ X\cong M_{k}(A) \textnormal{ for some } A\in\projk{k}{S^{\op}}\}$ and where $\projk{k}{\Sop}$ denotes the full subcategory of $\mod{\Sop}$ consisting of finitely presented $\Sop$-modules of projective dimension at most $k$.}
\erem

\subsection{$u$-cluster categories}

Cluster categories were introduced by Buan, Marsh, Reineke, Reiten and Todorov in \cite{BMRRT}. They were also introduced indepentently for the type $A$ case in \cite{Caldero}. The $u$-cluster category was first introduced by Bernhard Keller in \cite[Section 8.4]{Keller4}. Let $k$ be an algebraically closed field and $H$ be a finite dimensional hereditary $k$-algebra. For an integer $u\geqslant 1$, the \textit{$u$-cluster category} $\C$ is defined by $\D^{f}(H^{\op})/\tau^{-1}\Sigma^{u}$, where $\tau$ is the AR translation of $D^{f}(H^{\op})$ (see \cite{Assem} or \cite{Auslander}, for example) and $\Sigma$ is its suspension. Here $\D^{f}(H^{\op})$ is shorthand for $\D^{f}(\mod{H^{\op}})$.

By Keller, \cite[Section 4, Theorem]{Keller4}, the canonical projection functor $\pi:\D^{f}(H^{\op})\rightarrow\C$ is triangulated. Hence by composition with the inclusion functor we obtain a full embedding
\begin{equation}
\xymatrix{\mod{H^{\op}}\ar@{^{(}->}[r]^-{\iota}\ar@{^{(}->}[dr] & \D^{f}(H^{\op})\ar[d]^-{\pi} \\
 & \C,}
\label{embedding1}
\end{equation}
for $u\geqslant 2$.

Now let $H=kQ$ be the path algebra of a quiver $Q$ with no loops or oriented cycles. Then $H$ is an hereditary algebra, hence coherent, and $H\in\D^{f}(H^{\op})/\tau^{-1}\Sigma^{u}=\C$ is maximal $u$-orthogonal. In particular, we have:
\begin{itemize}
\item $\Hom{\C}{H}{\Sigma^{i} H}=0$ for $i=1,\ldots, u$;
\item $\Hom{\C}{H}{\Sigma^{-i} H}=0$ for $i=1,\ldots, u$.
\end{itemize}
In addition, we have that $\Hom{\C}{H}{H}\cong H$, therefore the endomorphism algebra has global dimension $1$, and for $u\geqslant 2$, $H$ satisfies the hypotheses of Theorem \ref{moorespectrafinite}. Therefore there exists a full embedding $M:\mod{H^{\op}}\rightarrow\M$ where $\M=\M_{2}$ in Definition \ref{myM}. Hence, composing with the inclusion functor, we have a full embedding
\begin{equation}
\xymatrix{\mod{H^{\op}}\ar@{^{(}->}[r]^-{M}\ar@{^{(}->}[dr] & \M\ar@{^{(}->}[d] \\
 & \C.}
\label{embedding2}
\end{equation}
 We claim that this embedding coincides with that of \eqref{embedding1} from \cite{Keller4}.
 
It is clear that embeddings \eqref{embedding1} and \eqref{embedding2} are canonically equivalent on $\mathsf{add}(H)$. Thus they are equivalent on $\proj{H^{\op}}$, so we only need to extend the equivalence to projective dimension $1$ since $H$ is hereditary. Let $A$ and $B$ be $H^{\op}$-modules of projective dimension $1$ and suppose we have a module homomorphism $a:A\rightarrow B$. Take projective resolutions of $A$ and $B$:
\begin{eqnarray*}
0\longrightarrow P_{1}\rightlabel{f_{1}} P_{0}\rightlabel{f_{0}} A\longrightarrow 0, \\
0\longrightarrow Q_{1}\rightlabel{g_{1}} Q_{0}\rightlabel{g_{0}} B\longrightarrow 0.
\end{eqnarray*}
By elementary homological algebra, see \cite{Hilton}, we can lift the homomorphism $a:A\rightarrow B$ to a commutative diagram:
\begin{equation}
\xymatrix{0\ar[r] & P_{1}\ar[r]^-{f_{1}}\ar@{-->}[d]_-{p_{1}} & P_{0}\ar[r]^-{f_{0}}\ar@{-->}[d]_-{p_{0}} & A\ar[r]\ar[d]_-{a} & 0 \\
0\ar[r] & Q_{1}\ar[r]_-{g_{1}} & Q_{0}\ar[r]_-{g_{0}} & B\ar[r] & 0.}
\label{cluster_lift}
\end{equation}
Since there is a natural isomorphism $\tau:M|_{\proj{H^{\op}}}\rightarrow\pi\circ\iota|_{\proj{H^{\op}}}$ there are commutative diagrams:
\begin{equation}
\xymatrix{MP_{1}\ar[r]^{Mf_{1}}\ar[d]_-{\tau_{P_{1}}}^-{\sim} & MP_{0}\ar[d]_-{\tau_{P_{0}}}^-{\sim} \\
\pi\circ\iota(P_{1})\ar[r]_-{\pi\circ\iota(f_{1})} & \pi\circ\iota(P_{0})}
\qquad \textnormal{and }
\xymatrix{MQ_{1}\ar[r]^{Mg_{1}}\ar[d]_-{\tau_{Q_{1}}}^-{\sim} & MQ_{0}\ar[d]_-{\tau_{Q_{0}}}^-{\sim} \\
\pi\circ\iota(Q_{1})\ar[r]_-{\pi\circ\iota(g_{1})} & \pi\circ\iota(Q_{0}).}
\label{cluster_naturalisomorphism}
\end{equation}
In addition, applying the functors $M$ and $\pi\circ\iota$ to diagram \eqref{cluster_lift} yields the following commutative diagrams, respectively:
\begin{equation}
\xymatrix{MP_{1}\ar[r]^-{Mf_{1}}\ar[d]_-{Mp_{1}} & MP_{0}\ar[r]^-{Mf_{0}}\ar[d]_-{Mp_{0}} & MA\ar[r]^-{h}\ar[d]_-{Ma} & \Sigma MP_{1}\ar[d]_-{\Sigma Mp_{1}} \\
MQ_{1}\ar[r]_-{Mg_{1}} & MQ_{0}\ar[r]_-{Mg_{0}} & MB\ar[r]_-{j} & \Sigma MQ_{1},}
\label{cluster_triangles1}
\end{equation}
and
\begin{equation}
\xymatrix{\pi\circ\iota(P_{1})\ar[r]^-{\pi\circ\iota(f_{1})}\ar[d]_-{\pi\circ\iota(p_{1})} & \pi\circ\iota(P_{0})\ar[r]^-{\pi\circ\iota(f_{0})}\ar[d]_-{\pi\circ\iota(p_{0})} & \pi\circ\iota(A)\ar[r]^-{\theta}\ar[d]_-{\pi\circ\iota(a)} & \Sigma (\pi\circ\iota(P_{1}))\ar[d]_-{\Sigma (\pi\circ\iota(p_{1}))} \\
\pi\circ\iota(Q_{1})\ar[r]_-{\pi\circ\iota(g_{1})} & \pi\circ\iota(Q_{0})\ar[r]_-{\pi\circ\iota(g_{0})} & \pi\circ\iota(B)\ar[r]_-{\phi} & \Sigma (\pi\circ\iota(Q_{1})).}
\label{cluster_triangles2}
\end{equation}

Note that diagram \eqref{cluster_triangles2} comes by virtue of Keller's theorem that the canonical projection functor is triangulated \cite[Section 4, Theorem]{Keller4}.

Applying the functor $M$ and $\pi\circ\iota$ to the top row of diagram \eqref{cluster_lift} and using the natural isomorphism highlighted in diagrams \eqref{cluster_naturalisomorphism} gives the following commutative diagram:
\begin{equation}
\xymatrix{MP_{1}\ar[r]^-{Mf_{1}}\ar[d]_-{\tau_{P_{1}}}^-{\sim} & MP_{0}\ar[r]^-{Mf_{0}}\ar[d]_-{\tau_{P_{0}}}^-{\sim} & MA\ar[r]^-{h}\ar@{-->}[d]_-{\sigma_{A}}^-{\sim} & \Sigma MP_{1}\ar[d]_-{\Sigma \tau_{P_{1}}}^-{\sim} \\
\pi\circ\iota(P_{1})\ar[r]_-{\pi\circ\iota(f_{1})} & \pi\circ\iota(P_{0})\ar[r]_-{\pi\circ\iota(f_{0})} & \pi\circ\iota(A)\ar[r]_-{\theta} & \Sigma (\pi\circ\iota(P_{1})),}
\label{cluster_triangles3}
\end{equation}
where $\sigma_{A}$ exists by an axiom of triangulated categories, see, for instance \cite[Definition 1.1.1 (TR3)]{Neeman2} and is an isomorphism by the Five Lemma for triangulated categories, see \cite[Proposition 1.1.20]{Neeman2}.

Likewise, one obtains the following diagram:
\begin{equation}
\xymatrix{MQ_{1}\ar[r]^-{Mg_{1}}\ar[d]_-{\tau_{Q_{1}}}^-{\sim} & MQ_{0}\ar[r]^-{Mg_{0}}\ar[d]_-{\tau_{Q_{0}}}^-{\sim} & MB\ar[r]^-{j}\ar@{-->}[d]_-{\sigma_{B}}^-{\sim} & \Sigma MQ_{1}\ar[d]_-{\Sigma \tau_{Q_{1}}}^-{\sim} \\
\pi\circ\iota(Q_{1})\ar[r]_-{\pi\circ\iota(g_{1})} & \pi\circ\iota(Q_{0})\ar[r]_-{\pi\circ\iota(g_{0})} & \pi\circ\iota(B)\ar[r]_-{\phi} & \Sigma (\pi\circ\iota(Q_{1})).}
\label{cluster_triangles4}
\end{equation}

Combining diagrams \eqref{cluster_triangles1}, \eqref{cluster_triangles2}, \eqref{cluster_triangles3}, and \eqref{cluster_triangles4} gives the following three-dimensional diagram:
\begin{equation}
\xymatrix@C-50pt@!C{ & \pi\circ\iota(P_{1})\ar[rr]\ar[dd]|\hole & & \pi\circ\iota(P_{0})\ar[rr]\ar[dd]|\hole & & \pi\circ\iota(A)\ar[rr]\ar[dd]|\hole & & \Sigma(\pi\circ\iota(P_{1}))\ar[dd] \\
MP_{1}\ar[rr]\ar[ur]^-{\sim}\ar[dd] & & MP_{0}\ar[rr]\ar[ur]^-{\sim}\ar[dd] & & MA\ar[rr]\ar@{-->}[ur]^-{\sim}\ar[dd] & & \Sigma MP_{1}\ar[ur]^-{\sim}\ar[dd] & \\
 & \pi\circ\iota(Q_{1})\ar[rr]|\hole & & \pi\circ\iota(Q_{0})\ar[rr]|\hole & & \pi\circ\iota(B)\ar[rr]|\hole & & \Sigma(\pi\circ\iota(Q_{1})) \\
MQ_{1}\ar[rr]\ar[ur]^-{\sim} & & MQ_{0}\ar[rr]\ar[ur]^-{\sim} & & MB\ar[rr]\ar@{-->}[ur]^-{\sim} & & \Sigma MQ_{1},\ar[ur]^-{\sim}}
\label{another_3d}
\end{equation}
where each square and cube is known to commutative except for those involving the broken arrows. We claim that the commutativity of the rest of the diagram forces the whole diagram to commute.

By the commutativity of the rest of the diagram, we have another diagram of distinguished triangles:
\begin{equation*}
\xymatrix{MP_{1}\ar[r]^-{Mf_{1}}\ar[d]_-{\gamma_{1}} & MP_{0}\ar[r]^-{Mf_{0}}\ar[d]_-{\gamma_{0}} & MA\ar[r]^-{h}\ar[d]_-{\gamma} & \Sigma MP_{1}\ar[d]_-{\Sigma \gamma_{1}} \\
\pi\circ\iota(Q_{1})\ar[r]_-{\pi\circ\iota(g_{1})} & \pi\circ\iota(Q_{0})\ar[r]_-{\pi\circ\iota(g_{0})} & \pi\circ\iota(B)\ar[r]_-{\phi} & \Sigma (\pi\circ\iota(Q_{1})),}
\label{cluster_triangles5}
\end{equation*}
where we have
\begin{eqnarray*}
\gamma_{1} & = & \pi\circ\iota(p_{1})\circ\tau_{P_{1}} - \tau_{Q_{1}}\circ Mp_{1} \\
\gamma_{0} & = & \pi\circ\iota(p_{0})\circ\tau_{P_{0}} - \tau_{Q_{0}}\circ Mp_{0} \\
\gamma 	   & = & \pi\circ\iota(a)\circ\sigma_{A} - \sigma_{B}\circ Ma.
\end{eqnarray*}
By the commutativity of diagrams \eqref{cluster_naturalisomorphism}, we have $\gamma_{1}=\gamma_{0}=\Sigma\gamma_{1}=0$. Thus we obtain the diagram:
\begin{equation*}
\xymatrix{MP_{1}\ar[r]^-{Mf_{1}} & MP_{0}\ar[r]^-{Mf_{0}}\ar[dr]_-{0} & MA\ar[r]^-{h}\ar[d]_-{\gamma} & \Sigma MP_{1}\ar@{-->}[dl]^-{\exists} \\
 &  & \pi\circ\iota(B). & }
\label{cluster_triangles6}
\end{equation*}
From diagram \eqref{cluster_triangles4} we have that $\pi\circ\iota(B)\cong MB$, thus $\pi\circ\iota(B)\in\M$. The broken arrow, $\Sigma MP_{1}\rightarrow \pi\circ\iota(B)$ must be zero by the proof of Theorem \ref{moorespectra}. Hence, it follows that $\gamma=0$, so that
$$\pi\circ\iota(a)\circ\sigma_{A} = \sigma_{B}\circ Ma.$$
This then forces diagram \eqref{another_3d} to commute, as claimed. In particular, we obtain the following commutative diagram for any module homomorphism $a:A\rightarrow B$ in $\mod{H^{\op}}$:
$$\xymatrix{MA\ar[r]^{Ma}\ar[d]_-{\sigma_{A}}^-{\sim} & MB\ar[d]_-{\sigma_{B}}^-{\sim} \\
\pi\circ\iota(A)\ar[r]_-{\pi\circ\iota(a)} & \pi\circ\iota(B).}$$
Hence there exists a natural isomorphism $\sigma: M\rightarrow \pi\circ\iota$ on $\mod{H^{\op}}$. We have, therefore, proved the following theorem.

\thm
Let $H=kQ$ be the path algebra of a quiver $Q$ with no loops or oriented cycles and let $\C=\D^{f}(H)/\tau^{-1}\Sigma^{u}$ be the $u$-cluster category as defined in \cite{Jorgensen_Holm} and \cite{Keller4} for $u\geqslant 2$. Then the canonical embedding, $\pi\circ\iota$, obtained in \cite[Section 4, Theorem]{Keller4} (see diagram \eqref{embedding1}) and the full embedding obtained in Theorem \ref{moorespectrafinite} (see diagram \eqref{embedding2}) are naturally equivalent.
\ethm

\

\noindent\textbf{Acknowledgement.} The author would like to thank his supervisor, Prof Peter J\o rgensen, for helpful discussions in the preparation of this paper and also an anonymous referee for suggesting a way to clarify the proof of Theorem \ref{moorespectra}.

\end{document}